\documentclass[11pt]{amsart}
\usepackage{amsfonts}
\usepackage{amsmath, amsthm, amssymb,color}
\usepackage{latexsym}

\numberwithin{equation}{section}
\allowdisplaybreaks

\usepackage{graphics}
\usepackage{graphicx}
\usepackage{epsfig}

\definecolor{darkblue}{rgb}{.1, 0.1,.8}
\definecolor{darkgreen}{rgb}{0,0.8,0.2}
\definecolor{darkred}{rgb}{.8, .1,.1}
\newcommand{\bfi}{\begin{fig}}
\newcommand{\efi}{\end{fig}}
\newcommand{\levy}{L\'evy}

\newtheorem{lemma}{Lemma}[section]

\newtheorem{theorem}[lemma]{Theorem}

\newtheorem{proposition}[lemma]{Proposition}
\newtheorem{definition}[lemma]{Definition}
\newtheorem{corollary}[lemma]{Corollary}
\newtheorem{example}[lemma]{Example}
\newtheorem{exercise}[lemma]{Exercise}
\newtheorem{remark}[lemma]{Remark}
\newtheorem{fig}[lemma]{Figure}
\newtheorem{tab}[lemma]{Table}

\newcommand{\bth}{\begin{theorem}}
\newcommand{\ethe}{\end{theorem}}

\newcommand{\bre}{\begin{remark}\em }
\newcommand{\ere}{\end{remark}}

\newcommand{\ble}{\begin{lemma}}
\newcommand{\ele}{\end{lemma}}

\newcommand{\pp}{point process}
\newcommand{\bde}{\begin{definition}}
\newcommand{\ede}{\end{definition}}
\newcommand{\bco}{\begin{corollary}}
\newcommand{\eco}{\end{corollary}}

\newcommand{\bpr}{\begin{proposition}}
\newcommand{\epr}{\end{proposition}}

\newcommand{\bexer}{\begin{exercise}}
\newcommand{\eexer}{\end{exercise}}

\newcommand{\bexam}{\begin{example}\rm }
\newcommand{\eexam}{\end{example}}

\newcommand{\btab}{\begin{tab}}
\newcommand{\etab}{\end{tab}}

\newcommand{\beao}{\begin{eqnarray*}}
\newcommand{\eeao}{\end{eqnarray*}\noindent}

\newcommand{\beam}{\begin{eqnarray}}
\newcommand{\eeam}{\end{eqnarray}\noindent}

\newcommand{\beqq}{\begin{equation}}
\newcommand{\eeqq}{\end{equation}\noindent}

\newcommand{\beali}{\begin{align*}}
\newcommand{\reeal}{\end{align*}}

\newcommand{\bce}{\begin{center}}
\newcommand{\ece}{\end{center}}

\newcommand{\barr}{\begin{array}}
\newcommand{\earr}{\end{array}}

\newcommand{\vague}{\stackrel{\lower0.2ex\hbox{$\scriptscriptstyle
                    \it{v} $}}{\rightarrow}}
\newcommand{\weak}{\stackrel{\lower0.2ex\hbox{$\scriptscriptstyle
                    \it{w} $}}{\rightarrow}}
\newcommand{\what}{\stackrel{\lower0.2ex\hbox{$\scriptscriptstyle
                    \it{\hat{w}} $}}{\rightarrow}}

\newcommand{\bdis}{\begin{displaymath}}
\newcommand{\edis}{\end{displaymath}\noindent}

\newcommand{\wt}{\widetilde}
\newcommand{\wh}{\widehat}

\newcommand{\bbf}{{\mathcal F}}

\newcommand{\ins}{insurance}

\renewcommand\d{{\mathrm d}}
\newcommand\e{{\mathrm e}}

\begin{document}

\bibliographystyle{alpha}
\title[Prediction in a non-homogeneous Poisson cluster model]
{Prediction in a non-homogeneous Poisson cluster model}
\today
\author[M. Matsui]{Muneya Matsui \vspace{3mm}\\ Nanzan University}
\address{Department of Business Administration, Nanzan University,
18 Yamazato-cho Showa-ku Nagoya, 466-8673, Japan}
\email{mmuneya@nanzan-u.ac.jp}

\begin{abstract}
 A non-homogeneous Poisson cluster model is studied, motivated by
 insurance applications. The Poisson center process 
 which expresses arrival times of claims, triggers off cluster member
 processes which correspond to number or amount of payments. The
 cluster member process is an additive process. Given the past observations of
 the process we consider expected values of future increments and
 their mean squared errors, aiming the application in claims reserving
 problems. 
 Our proposed process can cope with
 non-homogeneous observations such as the seasonality of claims arrival
 or the reducing property of payment processes, which are unavailable
 in the former models where both center and member processes are time
 homogeneous. 
 Hence results presented in this paper are
 significant extensions toward applications.  
 We also give numerical examples to show how non-homogeneity appears in predictions.
\end{abstract}

\vspace{.5ex}
\thanks{
This research was supported
by JSPS KAKENHI Grant number 25870879 Grant-in-Aid for Young Scientists (B)}

\pagecolor{white}

\maketitle 
\section{Introduction}
 A Cluster point process is one of the most important classes of point
 processes, which has two driving processes, the process of cluster
 center and the process of each cluster (see e.g. Daley and
 Vere-Jones \cite{daley:vere-jones:1988} or Westcott
 \cite{westcott:1971}). The Poisson cluster process is a version of
 cluster point processes whose process of center is a Poisson process. 
 The process has been applied to a wide-range of different fields such
 as earthquake aftershocks \cite{verejones:1970}, motor traffic
 \cite{bartlett:1963}, computer failure times \cite{lewis:1964} and
 broadband traffics \cite{hohn:veitch:abrey:2003} to name just
 a few. For more on the history and
 applications we refer to \cite{daley:vere-jones:1988}.

 Motivated by insurance applications we will investigate the
 Poisson cluster process of the form,
 \beam\label{eq:cpp}
 M(t)=\sum_{j=1}^{N(1)} L_{j}(t-T_j)\,, \quad t\ge 1\,,
 \eeam 
 where $0<T_1<T_2<\cdots$ are points of
 non-homogeneous Poisson ($NP$ for short) processes $N(t)$ and $(L_j),\,j=1,2,\ldots$ are 
 an iid sequence of additive processes with $L_j(t)=0,\,a.s.$ for
 $t\le 0$, such that $(T_j)$ and $(L_j)$ are independent. 
 In the insurance context $T_j \le 1$ would be the arrival of
 claims within a year, and $(L_{j}(t-T_j))_{j: T_j \le 1}$ are  
 the corresponding payment processes from an insurance company to
 policyholders. We could also regard the cluster as the counting process
 of payment number. Hence $M(t)$ would be the total number or amount of payments for the claims
 arriving in a year and being paid in the interval $[0,t],\, t\ge1$. Historically such kind of stochastic process modeling goes back to
 Lundberg (1903) (see comments in \cite[p.224]{mikosch:2009}) who introduced the Poisson process
 for a simple claim counting process. Norberg
 \cite{norberg:1993,norberg:1999} has been 
 considered to give publicity to the point process approach in
 a non-life insurance context.  
 
 Our focus in this paper is on the prediction of
 the future increments 
 \beao
 M(t,t+s]=M(t+s)-M(t)\,,\quad t\ge 1\,,s>0\,,
 \eeao
 for some suitable
 $\sigma$-fields $\bbf_t$ i.e. we will calculate $E[M(t,t+s]\mid
 \bbf_t]$ and evaluate the mean squared error of the prediction.  
 These kind of problems are known to be claims reserving problems, 
 which have been intensively studied from old times. 
 We refer e.g. to Chapter 11 of Mikosch \cite{mikosch:2009} for the recent
 development of the topic, where several interesting methods including famous {\it chain ladder} method
 are well explained. 

 In reference to prediction problems with the model
 \eqref{eq:cpp}, Mikosch \cite{mikosch:2009} introduced the model 
 into the claims reserving problems with a simple
 settings such that both the center process and clusters are
 homogeneous Poisson processes, where numerically 
 tractable form of predictor $E[M(t,t+s]\mid M(t)]$ is also obtained. More generally,
 Matsui and Mikosch \cite{matsui:mikosch:2010} consider L\'evy or truncated compound
 Poisson for clusters and obtain analytic forms of both prediction and its
 mean squared error. Matsui \cite{matsui:2011} introduced a variation of the model \eqref{eq:cpp}
 which starts randomly given number of cluster processes at each jump point of underlying
 process $N(t)$ and also obtain predictors and their errors. In a
 different path Jessen et al. \cite{jessen:mikosch:samorodnitsky:2009}
 takes simpler but useful point process modeling for the
 problem. See also Rolski and Tomanek
 \cite{rolski:tomanek:2011} which investigates asymptotics of conditional
 moments arising from prediction problems. Notice that almost all
 processes used in the the context are included in the class of L\'evy processes, which
 implies that increments are time-homogeneous. 
 
 In this paper we introduce non-homogeneity into both underlying Poisson
 process $N(t)$ and clusters $L_j$  by the use of 
 additive processes such that the processes have independent but
 not always stationary increments. More precisely, we assume
 a $NP$ 
 process for $N$, whereas each cluster $L_j$ is
 assumed to be an additive process which is given by a certain integral of
 a general Poisson random measure. 
 Our intention here is to model the seasonality of
 claims arrivals and the curved line of payment numbers or
 amounts which are naturally observed from data (see e.g.
 Table 2 of \cite{jessen:mikosch:samorodnitsky:2009}). 
 Again we emphasize that in the former models \cite{mikosch:2009}, \cite{matsui:mikosch:2010} or \cite{matsui:2011},
 they intensively use L\'evy clusters which are the processes of
 stationary independent increments and therefore are time
 homogeneous. 

 This paper is organized as follows.  
 In section \ref{sec:pipcm}, we consider the model with additive L\'evy processes and
 obtain the conditional characteristic function (ch.f. for abbreviation) $E[e^{ixM(t,t+s]}\mid M(t)]$. Based on
 the derived ch.f. we investigate expressions of $E[M(t,t+s]\mid M(t)]$
 where $NP$ 
 clusters and non-homogeneous negative
 binomial clusters are considered. In both cases, we derive recursive algorithm to
 calculate exact values of predictors and their conditional mean squared errors. In
 Section \ref{sec:predic:delay}
 the prediction $E[M(t,t+s]\mid \mathcal F_t]$ with different
 $\sigma$-fields $\mathcal F_t$ is investigated where we notice the
 delay in reporting times of claims and consider the number of reported
 claims until time $t$ for $\mathcal F_t$. Exact analytic forms for both
 predictors and their mean squared errors are calculated.
 In the final section, we give numerical examples to see how the 
 non-homogeneity affects the predictors. 

 Finally, we briefly explain basics of  
 an additive process $\{L(t)\}_{t\ge 0}$ based on Sato \cite[p.53]{sato:1999}. It is well known
 that the process is stochastically continuous, and has independent
 increments with c\`adl\`ag path starting at  
 $L(0)=0$ a.s. The distribution of the process $\{L(t)\}_{t\ge
 0}$ at time $t$ is determined by its generating triplet
 $(A_t,\nu_t,\gamma_t)$ since this determines the corresponding ch.f. 
 Among additive processes we work on the process of so called
 jump part such that the distribution of time $t$ is given by the
 inversion of 
 \beam
\label{chf:additive:process}
E[ \e^{ix L(t)}]=\exp\Big\{
\int_{(0,t]\times \mathbb{R}} ( \e^{ixv}-1)
\nu(\d (u,v))
\Big\},
\eeam
where 
a measure $\nu$ on $(0,\infty)\times \mathbb{R}$ satisfies
 $\nu((0,\cdot]\times \{0\})=0$, $\nu(\{t\}\times \mathbb{R})=0$ and
$$
\int_{(0,t]\times \mathbb{R}} (|u|\wedge 1) \nu(\d (u,v))<\infty,\qquad
 \mathrm{for}\ t \ge0.  
$$
In this case the
generating triplet is $(0,\nu_t,0)$ with $\nu_t(B):=\nu((0,t]\times B)$ for
any Borel set $B\in\mathcal B(\mathbb{R})$. 
The first condition means $\nu_t(\{0\})=0$ and the second one implies
stochastic continuity, whereas the third controls smoothness of the
path. 
In view of \eqref{chf:additive:process}
one see that an additive process has an integral representation by
 Poisson random measure on $(0,\infty)\times \mathbb{R}$ with intensity
 measure 
 $\nu$. 
We refer to Theorem 19.2 and
 19.3 of Sato \cite{sato:1999} for the jump part of an additive
 process. Although we could treat more general 
additive process by including the continuous part or another version of
 jump part, the
prediction procedure would be more complicated and we confine the 
process of the cluster as such.

\section{Prediction in Poisson cluster model}\label{sec:pipcm}
In this section firstly we give general prediction results which are valid for all 
additive L\'evy clusters given by \eqref{chf:additive:process} and then
we investigate numerically tractable expressions with examples. More
precisely, we study expressions of the conditional expectation of
$M(t,t+s]$ given $M(t)$, $t\ge 1,s>0$ and its
mean squared error. 

The main difference of our prediction from Matsui and Mikosch
\cite{matsui:mikosch:2010} is that we can not use the stationary
increments of cluster center nor cluster member processes and hence
 expressions for predictors require some devices and are more
 complicated. However, by discarding time 
 homogeneity of L\'evy processes, we can introduce time dependency into
 the process of the cluster, which is of critical importance in applications.  

The following is basic for the model \eqref{eq:cpp}. 
\ble
Assume the model \eqref{eq:cpp} with iid additive processes
$L_k,\,k=1,2,\ldots$ and a 
NP process $N$ with
mean measure $\Lambda$ such that $\Lambda[0,\infty)<\infty$. We write the generic of processes $L_k$ as
$L$. Then the ch.f. is given by 
\begin{align*}
E[\,\e^{ixM(t)}] = \exp\Big\{\int_{[0,1]} (E[\,\e^{ixL(t-u)}]-1)
 \Lambda(\d u) \Big\}
\end{align*}
for $t \ge 1$ and $x\in\mathbb{R}$. 
Moreover, assume that $E[L(t)]$ finitely exists for all $t \ge 0$, then
\begin{align*}
E[M(t)] =& \int_{[0,1]}E[L(t-u)]\Lambda(\d u),\quad t\ge 1.  
\end{align*}
Assume that $E[L^2(t)]$ is finite for all $t \ge 0$. Then, for $1\le s\le t$,  
\begin{align*}
\mathrm{Cov}(M(s),M(t)) =&
 \int_{[0,1]}(E[L^2(s-u)]+E[L(s-u)]E[L(s-u,t-u]])\Lambda(\d u). 
\end{align*}
\ele

We are starting to observe the conditional ch.f. of
$M(t,t+s]$ given $M(t)$. 

\ble \label{lem:characteristic:additive}
Assume the model \eqref{eq:cpp} with iid additive
processes $L_k,\,k=1,2,\ldots$ given by \eqref{chf:additive:process}
and a NP 
process $N$ with the mean value function
$\Lambda(\cdot)$. 
For $m=1,2,\ldots,s>0,\,t\ge 1$ and $x\in\mathbb{R}$, the conditional
ch.f. of $M(t,t+s]$ given $\{M(t)\in A\}$ for any Borel set $A$ has the
following form 
\beam \label{condi:chara:general}
&& \wh f_A(x) = E[\e^{ix M(t,t+s]} \mid M(t)\in A ]  \\
&=& \frac{
E\Big[ \exp\big\{
\sum_{j=1}^{N(1)} \int_{(t-V_j,t+s-V_j]\times \mathbb{R}}
(\e^{ixv}-1)\nu(\d (u,v))
\big\}
P\big( \sum_{j=1}^{N(1)}L_j(t-V_j)\in A  \mid
(V_j)
\big) \Big] }{
P\Big(\sum_{j=1}^{N(1)}L_j(t-V_j)\in A \Big)} \nonumber
\eeam
for an iid sequence $(V_j)$ with density function
\beam \label{eq:density-V}
F(\d x) &=& \Lambda(\d x)/\Lambda(1),\quad 0\le x\le 1 
\eeam
 such that $(V_j)$, $(L_j)$ and $N$ are mutually 
independent. 
\ele

\begin{proof}
Since $M(t)$ is measurable with respect to $\sigma$-filed by
 $(T_j),\,(L_j(t-T_j))$, we use the iteration property of conditional
 expectation to calculate 
\beao
&& E[\e^{ixM(t,t+s]}\mid M(t)] \\
&& = E\Big[
E\big[
\prod_{j=1}^\infty \e^{ixI_{\{T_j\le 1\}}L_j(t-T_j,t+s-T_j]} \mid (T_j), (L_j(t-T_j))
\big] \mid M(t)
\Big]\\
&& = E\Big[
E\big[
\prod_{j=1}^\infty \big( I_{\{T_j>1\}}+ I_{\{T_j\le 1\}} \e^{ix L_j(t-T_j,t+s-T_j]} \big)\mid (T_j), (L_j(t-T_j))
\big] \mid M(t)
\Big]\\
&& = E\Big[
\prod_{j=1}^{N(1)} E[\e^{ixL_j(t-T_j, t+s-T_j]} \mid T_j,L_j(t-T_j)] 
\mid M(t) \Big] \\
&& = E\Big[
\prod_{j=1}^{N(1)} E[\e^{ix L (t-T_j, t+s-T_j]} \mid T_j] 
\mid M(t) \Big]. 
\eeao
Accordingly, for any Borel set $A$, we obtain
\beao
&& E[\e^{ixM(t,t+s]}\mid M(t)\in A] \\
&& = \frac{
E\Big[ \prod_{j=1}^{N(1)} E[\e^{ixL (t-T_j, t+s-T_j]} \mid T_j] 
{\bf 1}_{\{ M(t)\in A \}} \Big] }{
P\Big( M(t) \in A \Big)} \\
&& = \frac{
E\Big[ E \big[ \prod_{j=1}^{N(1)} E[\e^{ixL (t-T_j, t+s-T_j]} \mid T_j] 
{\bf 1}_{\{ M(t) \in A \}} \mid
(T_j) \big] \Big] }{
P\Big(M(t)\in A \Big)} \\
&& =  \frac{
E\Big[ \prod_{j=1}^{N(1)} E[\e^{ixL (t-T_j, t+s-T_j]} \mid T_j] 
P\big(M(t) \in A  \mid
(T_j)
\big) \Big] }{
P\Big( M(t) \in A \Big)}.  
\eeao
Since quantities 
$$
\prod_{j=1}^{N(1)} E[\e^{ixL(t-T_j,t+s-T_j]}\mid T_j]\quad \mathrm{and}\quad
P\big(\sum_{j=1}^{N(1)} L_j(t-T_j) \in A \mid (T_j) \big)
$$
do not depend on the order of $(T_j)$, 
the order statistic property of Poisson yields
\beao
&& E\Big[ E\big[
\prod_{j=1}^{N(1)} E\big[\e^{ixL(t-T_j,t+s-T_j)}\mid T_j \big] P\big(
\sum_{j=1}^{N(1)} L_j(t-T_j)\in A \mid (T_j)
\big) \mid N(1) \big] \Big] \\
&& = E \Big[ \prod_{j=1}^{N(1)} E\big[\e^{ixL(t-V_j,t+s-V_j]}\mid V_j \big] 
P\big( \sum_{j=1}^{N(1)} L_j(t-V_j)\in A \mid (V_j)\big) \Big],
\eeao
where $(V_j)$ is the iid sequence whose common distribution is by \eqref{eq:density-V}. 
Now we insert this and \eqref{chf:additive:process} into the final expression
 and we obtain the result. 
\end{proof}
  
Based on $\hat f_A(x)$, we see important examples in the following
subsections.    

\subsection{Non-homogeneous Poisson clusters}\label{subsec:nhpc}
We consider the model of \eqref{eq:cpp} with $NP$ 
clusters $L_j,\,j=1,2,\ldots$ such that 
the generic cluster process $L$ at time $t$ 
has ch.f., 
\beam \label{eq:ch-poisson}
E[\e^{iuL(t)}]=\e^{\mu(t)(\e^{iu}-1)},\quad u\in \mathbb{R}. 
\eeam
where $\mu(t):=\nu((0,t]\times\{1\})$, i.e. the measure $\nu$ in \eqref{chf:additive:process}
has the support only on $(0,t]\times\{1\}$. The measure $\mu(t)$ is generally
called the mean value function or intensity measure of the Poisson
process 
(see Sec.19 of Sato \cite{sato:1999} or Sec.7.2 of Mikosch
\cite{mikosch:2009}). By the condition $\nu(\{t\}\times \mathbb{R})=0$ 
(stochastic continuity) before, $\mu(t)$ is assumed to be continuous in
$t$. 
Moreover, we assume that $\mu(0,\infty)<\infty$. Notice that  
the Poisson
process is one of the most important processes among
additive processes. Besides, it is a basic process for modeling the claim 
reserves
in the non-life insurance context. 
 
Before constructing prediction we define some notations. Let
$\phi_{(Y_1,Y_2)}(z_1,z_2)$ be the Laplace transform of a bivariate
random variable $(Y_1,Y_2),\,
\phi_{(Y_1,Y_2)}(z_1,z_2):=E[\e^{-z_1Y_1-z_2Y_2}],\ z_1\ge0, z_2\ge0$ and
denote its $(m,n)$th
partial derivatives by $\phi^{(m,n)}_{(Y_1,Y_2)}(z_1,z_2)$, 
whereas $\phi^{(m)}_Y(t),\ldots$ denotes
simply $m$th derivative of $\phi_Y(z)=E[\e^{-zY}]$ with $z\ge0 $.  
Throughout we use a random sum
\beam \label{sum:mu}
R_{N(1)}(t):= \sum_{j=1}^{N(1)} \mu(t-V_j),\quad t\ge1, 
\eeam
where 
$(V_j)$ is an iid random
sequence with common density \eqref{eq:density-V}. 

\ble
Assume the model \eqref{eq:cpp} with iid $NP$ 
processes $L_k,\,k=1,2,\ldots$ with mean value function $\mu(\cdot)$ and a
NP process $N$ with the intensity measure
$\Lambda(\cdot)$. Then for $m=1,2,\ldots$ and $x\in\mathbb{R}$ the conditional
ch.f. of $M(t,t+s],\,t\ge1,s>0$ given $\{M(t)=m\}$ has the
following form 
\beam \label{condi:chara:nonhomo}
\wh f_m(x) &=& E[\e^{ix M(t,t+s]} \mid M(t)=m ] \\
&=& \frac{ \phi_{R_{N(1)}(t), R_{N(1)}(t+s)}^{(m,0)}(\e^{ix},1-\e^{ix})}{\phi_{R_{N(1)}}^{(m)}(1)},
\nonumber
\eeam
where the random element $R_{N(1)}$ is given by \eqref{sum:mu}. 
\ele

\begin{proof}
By inserting the ch.f.
 \eqref{eq:ch-poisson} into the
 expression \eqref{condi:chara:general} in Lemma \ref{lem:characteristic:additive}, we observe 
\beao
&& E[\e^{ixM(t,t+s]}\mid M(t)=m] \\
&& =  \frac{
E\Big[ \e^{\sum_{j=1}^{N(1)}\mu(t-T_j,t+s-T_j](\e^{ix}-1)}
P\big( \sum_{j=1}^{N(1)} L_j(t-V_j)=m  \mid (V_j)
\big) \Big] }{
P\Big(\sum_{j=1}^{N(1)}L_j(t-V_j)=m \Big)}. 
\eeao 
The aggregation property of Poisson processes (Prop. 7.3.11 of \cite{mikosch:2009}) yields
\beao
&& P\Big( \sum_{j=1}^{N(1)}L_j(t-V_j)=m  \mid
(V_j) \Big) = 
\frac{(
\sum_{j=1}^{N(1)}\mu(t-V_j) )^m }{m!} \e^{-\sum_{j=1}^{N(1)}\mu(t-V_j) }\
 a.s.
\eeao
from which it follows that
\beao
&& \e^{\sum_{j=1}^{N(1)}\mu(t-V_j,t+s-V_j]( \e^{ix}-1)}
P\big( \sum_{j=1}^{N(1)} L_j(t-V_j)=m  \mid (V_j)
\big) \\
&& = \frac{(\sum_{j=1}^{N(1)}\mu(t-V_j))^m}{m!} \e^{ \sum_{j=1}^{N(1)}\{
 \mu(t+s-V_j)(\e^{ix}-1)
-\mu(t-V_j)\e^{ix}\}
} \\
&& = \frac{(-1)^m}{m!} \big(
 \e^{ \sum_{j=1}^{N(1)}\{
 \mu(t+s-V_j)(\e^{ix}-1)
- \mu(t-V_j)y\}
}
\big)_y^{(m)} \mid_{y=\e^{ix}}. 
\eeao
Now taking expectation for $(V_j)$, we obtain by Fubini's theorem that 
\begin{align*}
& P\Big(\sum_{j=1}^{N(1)}L_j(t-T_j) = m \Big) = \frac{(-1)^m
 \phi^{(m)}_{R_{N(1)}(t)}(1)}{m!}, \\
& E\Big[ \e^{ \sum_{j=1}^{N(1)}\mu(t-V_j, t+s-V_j](\e^{ix}-1)}
 P\big( \sum_{j=1}^{N(1)}L_j(t-T_j) = m  \mid
(V_j)
\big) \Big] \\
& = \frac{(-1)^m}{m !}
\big( \phi_{R_{N(1)}(t), R_{N(1)}(t+s)}(y,1-\e^{ix}) \big)_y^{(m)} \mid_{y=\e^{ix}}
\end{align*}
for $x\in \mathbb{R}$.
Hence we obtain the result. 
\end{proof}

Now by differentiating \eqref{condi:chara:nonhomo} sufficiently often, we obtain the conditional moments. 

\bth
Assume the model \eqref{eq:cpp} with iid $NP$ 
processes
$L_k,\,k=1,2,\ldots$ with the mean value function $\mu(\cdot)$ and a
$NP$ process $N$ with the mean value
function $\Lambda(\cdot)$. Then the prediction $\widehat M_m(t,t+s]$ of
$M(t,t+s]$ given $\{M(t)=m\}$ has the form 
\beam\label{cond:mean1}
\widehat M_m(t,t+s] =
\frac{\phi^{(m+1,0)}_{R_{N(1)}(t),R_{N(1)}(t+s)}(1,0)
- \phi^{(m,1)}_{R_{N(1)}(t),R_{N(1)}(t+s)}(1,0)
}{\phi^{(m)}_{R_{N(1)}(t)}(1)}
\eeam
and the conditional variance of $M(t,t+s]$ given $\{M(t)=m\}$ is 
\beam\label{cond:vari1}
&& \mathrm{Var}(M(t,t+s] \mid M(t)=m) \\
&&= 
\frac{\phi^{(m+2,0)}_{R_{N(1)}(t),R_{N(1)}(t+s)}(1,0) - 2 \phi^{(m+1,1)}_{R_{N(1)}(t),R_{N(1)}(t+s)}(1,0) 
+ \phi^{(m,2)}_{R_{N(1)}(t),R_{N(1)}(t+s)}(1,0)
}{\phi^{(m)}_{R_{N(1)}(t)}(1)} \nonumber \\
&&\qquad +\widehat M_m(t,t+s] -( \widehat M_m(t,t+s])^2
, \nonumber
\eeam
where 
$R_{N(1)}$ is the random sum \eqref{sum:mu}. 
\ethe

It is desirable to obtain an explicit expression for the unconditional
mean squired error $E[(M(t,t+s]-E[M(t,t+s]\mid M(t)])^2]$, 
since it gives a certain measure for evaluating goodness of
predictors. However, in the light of expressions \eqref{cond:mean1} and \eqref{cond:vari1} it seems
intractable (see Remark 2.2 of Matsui and Mikosch
\cite{matsui:mikosch:2010}.). 
Hence we content with conditional moments which are provided with numerically tractable expressions. 

In what follows we investigate further expressions of 
\eqref{cond:mean1} and \eqref{cond:vari1}. It is convenient to observe the bivariate Laplace
transform of $(R_{N(1)}(t), R_{N(1)}(t+s))$ i.e.
\beam \label{ch-f:bivariate}
\phi_{R_{N(1)}(t), R_{N(1)}(t+s)}(y,z) &=& E[
\e^{-y R_{N(1)}(t)-z R_{N(1)}(t+s) }] \\
&=& E\Big[
\prod_{j=1}^{N(1)}
E[
\e^{-\mu(t-V_j)y-\mu(t+s-V_j)z} ] \mid N(1)
\Big] \nonumber \\
&=& E\Big[
\Big(
\int_0^1 \e^{-\mu(t-v)y-\mu(t+s-v)z}  \frac{\Lambda(\d v)}{\Lambda (1)}
\Big)^{N(1)}
\Big] \nonumber \\
&=&
\exp \Big\{
\int_0^1 ( \e^{-\mu(t-v)y-\mu(t+s-v)z}-1) \Lambda(\d v) \Big\} \nonumber \\
&=&
\e^{\Lambda(1)(\phi_{R_1(t),R_1(t+s)}(y,z)-1)} \nonumber
\eeam
and derivatives of
$\Lambda(1)\phi_{R_1(t),R_1(t+s)}(y,z)$ with respect to $z$ at $z=0$, 
$$
\psi_j(y)=\Lambda(1)\phi^{(0,j)}_{R_1(t),R_1(t+s)}(y,0)=\int_0^1(-\mu(t+s-v))^j
\e^{-\mu(t-v)y}\Lambda(\d v),\quad j=0,1,2,
$$
where we note that $\Lambda(1)\phi_{R_{N(1)}(t)}(y)=\psi_0(y)$.

\ble
Let $\ell=1,2,\ldots$ and $j=0,1,2$. 
Let $\phi^{(\ell,j)}_{R_{N(1)}(t),R_{N(1)}(t+s)}(1,0)$ be the $(\ell,j)$th
partial derivative of $\phi_{R_{N(1)}(t),R_{N(1)}(t+s)}(y,z)$ at $(y,z)=(1,0)$ and let
$\psi^{(\ell)}_j(1)$ be the $\ell$th derivative of $\psi_j(y)$ at
$y=1$. Then, the following recursive relations hold. 
\beao
\psi_j^{(\ell)}(1) = \int_0^1 (-\mu(t+s-v))^j(-\mu(t-v))^\ell
\e^{-\mu(t-v)}\Lambda(\d v),\quad j=0,1,2, 
\eeao
and 
\beao
\phi^{(\ell,0)}_{R_{N(1)}(t),R_{N(1)}(t+s)}(1,0) &=& \sum_{k=0}^{\ell-1}
\binom{\ell-1}{k}
\phi_{R_{N(1)}(t),R_{N(1)}(t+s)}^{(k,0)}(1,0)\psi_0^{(\ell-k)}(1),\\
\phi^{(\ell,1)}_{R_{N(1)}(t),R_{N(1)}(t+s)}(1,0) &=& \sum_{k=0}^{\ell}
\binom{\ell}{k}
\phi_{R_{N(1)}(t),R_{N(1)}(t+s)}^{(k,0)}(1,0)\psi_1^{(\ell-k)}(1),\\
\phi^{(\ell,2)}_{R_{N(1)}(t),R_{N(1)}(t+s)}(1,0) &=&
\sum_{k=0}^{\ell}
\binom{\ell}{k}
\phi_{R_{N(1)}(t),R_{N(1)}(t+s)}^{(k,1)}(1,0)
\psi_1^{(\ell-k)}(1) \\
&& + \sum_{k=0}^{\ell}
\binom{\ell}{k} \phi_{R_{N(1)}(t),R_{N(1)}(t+s)}^{(k,0)}(1,0)\psi_2^{(\ell-k)}(1).
\eeao
\ele

\begin{proof}
We differentiate ch.f. \eqref{ch-f:bivariate} to see
\beao
\phi_{R_{N(1)}(t),R_{N(1)}(t+s)}^{(1,0)}(y,0) &=& \phi_{
 R_{N(1)}(t),R_{N(1)}(t+s)}(y,0)\psi_0^{(1)}(y), \\
\phi_{R_{N(1)}(t),R_{N(1)}(t+s)}^{(0,1)}(y,0) &=& \phi_{
 R_{N(1)}(t),R_{N(1)}(t+s)}(y,0)\psi_1(y), \\
\phi_{R_{N(1)}(t),R_{N(1)}(t+s)}^{(0,2)}(y,0) &=& 
\phi^{(0,1)}_{
 R_{N(1)}(t),R_{N(1)}(t+s)}(y,0) \psi_1(y) \\
&& + 
\phi_{
 R_{N(1)}(t),R_{N(1)}(t+s)}(y,0)\psi_2(y). \\
\eeao
Applications of the Leibniz's rule to these quantities yield our desired
 results. 
\end{proof}

\subsection{Non-homogeneous negative binomial clusters}
We consider a negative binomial ($NB$ for short) process for the generic random process $L$
of clusters $L_k,\,k=1,2,\ldots$ 
Let $\mu(t)>0$ be a continuous function of time $t$ with $\mu(0)=0$ and
let $p\in(0,1)$ so that $q=1-p$. 
\eqref{chf:additive:process} is constructed with $\nu(\d (u,v))=\mu(\d
u) \times \sigma(\d v)$ where $\sigma$ is concentrated on positive
integer and $\sigma(\{k\})=k^{-1}q^k,\,k=1,2,\ldots$ Accordingly the
ch.f. of $L(t)$ has the form 
\begin{align}
\label{ch:NB}
E[\e^{ix L(t)}]=\left(\frac{p}{1-q \e^{ix}}\right)^{\mu(t)},\quad x\in \mathbb{R}
\end{align}
(see e.g. \cite[p.20]{sato:1999}). The marginal distribution at time $t$ of the
process follows $NB$ with parameters $\mu(t)$ and $p$
(we also write $NB(\mu(t),p)$ for abbreviation) i.e.
\begin{align}
\label{pb:NB}
P(L(t)=k)=\binom{\mu(t)+k-1}{k} p^{\mu(t)}
q^k \ k=1,2,\ldots
\end{align}
such that the mean and variance of the process are respectively given by 
$E[L(t)]=\mu(t)q/p$ and $\mathrm{Var}(L(t)) =\mu(t) q/p^2$.
The distributions of increments $L(t)-L(s)$, $0\le s\le t<\infty$ are 
mutually independent and follows $NB(\mu(t)-\mu(s),p)$. Although there exist only a few
references for non-homogeneous $NB$ process, e.g. Carrillo
\cite{carrilo:1989}, 
for homogeneous $NB$ process, 
detailed distributional properties
 are given in e.g. Kozubowski and Podgorski \cite{kozubowski:podgorski:2009} 
(see also Johnson et
al. \cite{johnson:kotz:kemp:1992}). 

Throughout we use the bivariate probability generating function 
$G_{Y_1,Y_2}(z_1,z_2):=E[z_1^{Y_1}z_2^{Y_2}],\,|z_1z_2|\le 1,\,
(z_1,z_2)\in \mathbb{C}^2$
of an
integer valued random vector $(Y_1,Y_2)$, 
and its $(k,\ell)$th derivatives
$G^{(k,\ell)}_{Y_1,Y_2}(z_1,z_2)$ with respect to $(z_1,z_2)$, $\{k,\ell\}=1,2,\ldots$ 
Moreover, the notation $(\cdot)_z^{(m)}$ denotes the $m$th derivative of the
quantity in the brace. We again use the random sum \eqref{sum:mu} where
$\mu$ is replaced by that of $NB(\mu(t),p)$. 
We abbreviate $G_{R_{N(1)}(t),R_{N(1)}(t+s) }(z_1,z_2)$ to
$G_{t,t+s}(z_1,z_2)$ throughout this section.

It is convenient to start with the conditional ch.f. of 
$M(t,t+s]$ given $M(t)$. 

\ble
Assume the model \eqref{eq:cpp} with iid $NB(\mu(t),p)$ additive clusters
$L_k,\,k=1,2,\ldots$ such that $\mu(t)>0$ is continuous 
and $p\in (0,1)$. 
Then conditional ch.f. of $M(t,t+s]$ given
$\{N(t)=m\}$ has the form, 
\begin{align*}
\wh f_m(x)= \frac{1}{p(p^{m-1}G_{R_{N(1)}}(p))_p^{(m)}} \sum_{k=0}^{m-1} 
\binom{m}{k}\frac{(m-1)!}{(m-k-1)!}
(1-q\e^{ix})^{m-k} G_{t,t+s}^{(m-k,0)}(1-q\e^{ix},p/(1-q\e^{ix})), 
\end{align*}
where $R_{N(1)}$ is the random sum of \eqref{sum:mu} and $q=1-p$.
\ele

\begin{proof}
We apply \eqref{ch:NB} to the expression 
\eqref{condi:chara:general} of
 Lemma 2.2, which yields the
 conditional ch.f. for $NB$ clusters as
\begin{align*}
& E[\e^{ixM(t,t+s]}\mid M(t)=m] \\
& =  \frac{
E\Big[ (p/(1-q \e^{ix}))^{\sum_{j=1}^{N(1)}\mu(t-V_j,t+s-V_j]}
P\big( \sum_{j=1}^{N(1)}L_j(t-V_j) = m  \mid
(V_j)
\big) \Big] }{
P\Big(\sum_{j=1}^{N(1)}L_j(t-V_j) = m \Big)}.
\end{align*}
Since $N(1)$ is measurable with the $\sigma$-filed by $(V_j)$ and
since $(L_j)$ is independent of $(V_j)$, it follows from \eqref{pb:NB}  that 
\begin{align*}
 P\big(\sum_{j=1}^{N(1)}L_j(t-V_j) = m \mid (V_j) \big) &= 
\frac{p(1-p)^m}{m!} \frac{(R_{N(1)}+m-1)!}{(R_{N(1)}-1)!}p^{R_{N(1)}-1}
  \\
& = \frac{p(1-p)^m}{m!}(p^{R_{N(1)}+m-1})^{(m)}_p.
\end{align*}
We apply a similar calculation to the enumerator to obtain
\begin{align*}
& (p/(1-q\e^{ix}))^{\sum_{j=1}^{N(1)}\mu(t-V_j,t+s-V_j]}
P\big( \sum_{j=1}^{N(1)}L_j(t-V_j)=m  \mid
(V_j) \big) \\
& = (\gamma^{R_{N(1)}(t)+m-1})^{(m)}_\gamma \mid_{\gamma= 1-q\e^{ix}}
(1-q\e^{ix})  (p/(1-q\e^{ix}))^{R_{N(1)}(t+s)} \frac{(1-p)^m}{m!} \\
& = \frac{(1-p)^m}{m!}(1-qe^{ix})\Big(
\gamma^{m-1}G_{t,t+s}\big(\gamma,p/(1-q\e^{ix})\big)
\Big)_{\gamma}^{(m)} \mid_{\gamma=1-q\e^{ix}}.
\end{align*}\noindent
Now taking expectation for both quantities under notations of differentiation, which is
 justified by Fubini's theorem, we conclude the result. 
\end{proof}

Differentiation of the conditional ch.f. at $x=0$
several times yields the following result. 

\bth
Let $L$ be $NB(\mu(t),p)$ process such that $\mu(t)>0$ is continuous and
$p\in (0,1)$. 
Then the prediction $M(t,t+s]$ given
$\{M(t)=m\}$ is 
\begin{align*}
\wh M_m(t,t+s] &= \frac{1}{p (p^{m-1}G_{R_{N(1)}}(p))^{(m)}_p } 
\sum_{k=0}^{m-1} \binom{m}{k}\frac{(m-1)!}{(m-k-1)!} (-p^{m-k-1}q)\\
&\hspace{3cm} \times 
\Big\{
(m-k)G_{t,t+s}^{(m-k,0)}(p,1) +p G_{t,t+s}^{(m-k+1,0)}(p,1)-G_{t,t+s}^{(m-k,1)}(p,1)
\Big\}
\end{align*}
and the conditional mean squared error has the form 
\begin{align*}
& \mathrm{Var}(M(t,t+s]\mid M(t)=m) \\
& = -\wh M^2_m(t,t+s] + \frac{1}{p(p^{m-1}G_{R_{N(1)}}(p))^{(m)}_p } 
\sum_{k=0}^{m-1} \binom{m}{k}\frac{(m-1)!}{(m-k-1)!} p^{m-k-2}q\\
& \hspace{4cm} \times \Big[
(m-k)\{(m-k)q-1\} G_{t,t+s}^{(m-k,0)}(p,1) \\
&\hspace{4.5cm} + p\{2(m-k)q-p\}
G_{t,t+s}^{(m-k+1,0)}(p,1) \\
&\hspace{4.5cm} -\{2(m-k)q-1-q\}G_{t,t+s}^{(m-k,1)}(p,1) \\
& \hspace{4.5cm} + q \big\{ p^2 G_{t,t+s}^{(m-k+2,0)}(p,1) -2p
G_{t,t+s}^{(m-k+1,1)}(p,1)+ G_{t,t+s}^{(m-k,2)}(p,1)\big\}
\Big]. 
\end{align*}
\ethe

For numerical purpose, it is desirable to obtain tractable forms of
$G_{t,t+s}^{(k,j)}(p,1),\,j=0,1,2,\,k=1,2,\ldots$ Recall 
that $G_{t,t+s}(z_1,z_2)=G_{R_{N(1)}(t), R_{N(1)}(t+s)}(z_1,z_2)$ and one
see easily that 
\beao
G_{t,t+s}(z_1,z_2) &=& E[z_1^{R_{N(1)}(t)}z_2^{R_{N(1)}(t+s)}] \\
&=& E\big[
\big(
E[z_1^{\mu(t-V)}z_2^{\mu(t+s-V)}]
\big)^{N(1)}
\big] \\
&=& E\Big[
\Big(
\int_0^1 z_1^{\mu(t-v)}z_2^{\mu(t+s-v)}\frac{\Lambda(\d v)}{\Lambda(1)}
\Big)^{N(1)}
\Big] \\
&=& \e^{\Lambda(1)(G_{R_1(t),R_1(t+s)}(z_1,z_2)-1 )} 
\eeao
and derivatives of $\Lambda(1)G_{R_1(t),R_1(t+s)}(z_1,z_2)$ with
respect to $z_2$ at $z_2=1$,
\beam \label{derivative-H}
H_j(z_1):=\Lambda(1)G_{R_1(t),R_1(t+s)}^{(0,j)}(z_1,1) =\int_0^1
\frac{\Gamma(\mu(t+s-v)+1)}{\Gamma(\mu(t+s-v)+1-j)}
z_1^{\mu(t-v)} \Lambda(\d v),
\eeam
where $\Gamma(\cdot)$ is the Gamma function. Then we have

\bpr
Let $k,\ell=1,2,\ldots$ and $G^{(k,\ell)}_{t,t+s}(p,1)$ be the
$(k,\ell)$th derivatives of $G_{t,t+s}(z_1,z_2)$ at
$(z_1,z_2)=(p,1)$ and let $H_j^{(k)}(p),\,j=0,1,2$ be the $k$th
derivative of $H_j(z_1)$ at $z_1=p$. Then the following recursive relations hold.
\[
H_j^{(\ell)}(p) = \int_0^1 
\frac{\Gamma(\mu(t+s-v)+1)}{\Gamma(\mu(t+s-v)+1-j)} 
\frac{\Gamma(\mu(t-v)+1)}{\Gamma(\mu(t-v)+1-\ell)}
 z_1^{\mu(t-v)} \Lambda(\d v)
\]
and
\begin{align*}
G_{t,t+s}^{(\ell,0)}(p,1) &= \sum_{k=0}^{\ell-1} \binom{\ell-1}{k}
G_{R_{N(1)}(t)}^{(k)}(p)\cdot H_0^{(\ell-k)}(p) \\
G_{t,t+s}^{(\ell,1)}(p,1) &= \sum_{k=0}^\ell \binom{\ell}{k}
G_{R_{N(1)}(t)}^{(k)}(p)\cdot H_1^{(\ell-k)}(p) \\
G_{t,t+s}^{(\ell,2)}(p,1) &= \sum_{k=0}^\ell \binom{\ell}{k}
\Big\{ 
G_{t,t+s}^{(k,1)}(p)\cdot 
H_1^{(\ell-k)}(p)
+G_{R_{N(1)}(t)}^{(k)}(p) \cdot
H_2^{(\ell-k)}(p) \Big\}.
\end{align*}
\epr

\begin{proof}
We differentiate $G_{t,t+s}(z_1,z_2)$ with respect to $z_1$ and $z_2$
 proper times at $z_1=1$ and obtain 
\beao
G_{t,t+s}^{(1,0)}(z_1,1) &=& 
G_{R_{N(1)}(t)} (z_1) \cdot H_0^{(1)}(z_1) \\
G_{t,t+s}^{(0,1)}(z_1,1) &=& 
G_{R_{N(1)}(t)}(z_1)\cdot H_1(z_1) \\
G_{t,t+s}^{(0,2)}(z_1,1) &=& 
G^{(0,1)}_{t,t+s}(z_1,1)\cdot H_1(z_1)+  G_{R_{N(1)}(t)}(z_1)\cdot H_2(z_1).
\eeao
Applications of Leibniz's rule to these quantities together with \eqref{derivative-H}
yield the result. 
\end{proof}

\section{Prediction with delay in reporting}\label{sec:predic:delay}
In this section we introduce the time difference $D_k > 0$ between the
arrival time $T_k$ of $k$th claim and its reporting time, i.e.
the report of $k$th claim is coming at time $T_k+D_k$, and then we 
start the cluster process $L_k$. Accordingly in the model of \eqref{eq:cpp},
$L(t-T_k)$ are replaced by $L(t-(T_k+D_k))$ and we will work with model
\[
M(t)=\sum_{j=1}^\infty I_{\{T_j\le 1\}} L_j(t-(T_j+D_j)),\quad T_j\ge
1,\,D_j\ge0 . 
\]
We assume that
the generic random element $D$ of iid sequences $(D_k)$ takes positive values with
common distribution $F_D$ such that $(D_k)$ is independent of $(L_k)$ and $N$. 

Recall that usually the total claims number $N(1)$ may not be available at time $t\ge1$, while we
know the reported number of claims,
\begin{align*}
\wh N(t)=\# \{k\ge 1: T_k+D_k \le t, T_k \in [0,1] \}. 
\end{align*}
In what follows, we will consider the prediction $M(t,t+s]$ based on $\wh
N(t)$, namely we will calculate the conditional expectation
\[
\wh M_\ell (t,t+s] = E[M(t,t+s]\mid \wh N(t)=\ell],\, \ell=0,1,2,\ldots
\]
First we specify the distribution of $\wh N(t)$. Let $Q$ be a Poisson
random measure on the space $E=[0,1]\times [0,\infty)$ with mean measure
$\nu=\Lambda\times F_D$. Then $N(1)$ and $\wh N(t)$ have the Poisson
integral representation (c.f. Ex. 7.3.6 in Mikosch \cite{mikosch:2009}),
\begin{align}
\label{eq:pinteg:ovserved:number}
N(1) &= \int_E Q(\d s,\d y) 
 = \int_{s=0}^1 \int_{y=0}^{t-s} Q(\d s,\d y)+\int_{s=0}^1
 \int_{y=t-s}^\infty Q(\d s,\d y)\\
 &= \wh N(t) - [N(1)-\wh N(t)],
\end{align}
where random variables $N(1)-\wh N(t)$ and $\wh N(t)$ are independent
and Poisson distributed with parameters 
\[
E[N(1)-\wh N(t)] = \int_0^1 \int_{t-s}^\infty \Lambda(\d s)F(\d y) =\wh
\Lambda(t)\qquad \text{and}\qquad E[\wh N(t)]=\Lambda(1)-\wh \Lambda(t).
\]

It is convenient to start with the conditional ch.f. 
$M(t,t+s]$ given $\wh N(t)$. 

\ble
Let $\ell=0,1,2,\ldots,\,t\ge 1$ and $s>0$. The conditional
ch.f. of $M(t,t+s]$ given $\{\wh N(t)=\ell \}$ has the
form
\begin{align}
E[\e^{ixM(t,t+s]}\mid \wh N(t)=\ell] &=
(E[\e^{ixL(t-Z,t+s-Z)}])^{\ell} \nonumber \\
& \times \exp\Big\{
-\int_{v=0}^1 \int_{r=t-v}^{t+s-v} (1-E[
\e^{ixL(t+s-v-r)}
]) \Lambda(\d v) F_D(\d r)
\Big\},\label{condi:ch:poisson:observed}
\end{align}
where $Z$ has distribution $\Lambda*F_D/E[\wh N(t)]$. 
\ele

\begin{proof}
Since $\wh N(t)$ is measurable with respect to $\sigma$-field by $(T_j)$
 and $(D_j)$, the conditional ch.f. of
 $M(t,t+s]$ on $(T_j)$ and $(D_j)$ has the form,
\begin{align*}
& E[\e^{ixM(t,t+s]}\mid (T_j), (D_j)] \\
& = E\Big[
\prod_{k=1}^\infty \exp\{ixL_k(t-T_k-D_k,t+s-T_k-D_k]1_{\{T_k\le 1\}}\}
 \mid (T_j),(D_j)
\Big] \\
& = \prod_{k:T_k \le 1, T_k+D_k \le t}^\infty E[\e^{ix L(t-T_k-D_k,
 t+s-T_k-D_k]}\mid T_k,D_k] \\
& \quad \times \prod_{k:T_k \le 1, T_k+D_k > t}^\infty E[\e^{ix L(t+s-T_k-D_k)}\mid T_k,D_k], 
\end{align*}
where in the last step we notice $L(t-T_k-D_k)=0\ a.s.$ for
 $k:\,T_k+D_k\ge t.$ We proceed calculation by the chain rule of conditional
 expectation to obtain 
\begin{align*}
& E[\e^{ixM(t,t+s]}\mid \widehat N(t)] \\
& = E\Big[
\exp\Big\{
\sum_{k:T_k \le 1, T_k+D_k \le t} \log E[
\e^{ixL(t-T_k-D_k, t+s-T_k-D_k]} \mid T_k,D_k
] 
\Big\} \\
& \qquad \times
\exp\Big\{
\sum_{k:T_k \le 1, T_k+D_k > t} \log E[
\e^{ixL(t+s-T_k-D_k)} \mid T_k,D_k
] \Big\} 
 \mid \wh N(t)
\Big] \\
& = E\Big[
\exp\Big\{
\int_{v=0}^1 \int_{r=0}^{t-v} \log E[
\e^{ixL(t-v-r, t+s-v-r]}] Q(\d v,\d r)
\Big\} \\
& \qquad \times
\exp\Big\{
\int_{v=0}^1 \int_{r=t-v}^{t+s-v} \log E[
\e^{ixL(t+s-v-r)}] Q(\d v,\d r) \Big\}
 \mid \wh N(t)
\Big]. 
\end{align*}
In the last expression, the Poisson integral of $\wh N(t)$ and the
 second integral have disjoint support and hence they are
 independent. This together with the order statistics property of
 points of 
 $\wh N(t)$ yields 
\begin{align*}
& E[\e^{ixM(t,t+s]}\mid \wh N(t)] \\
& = E\Big[
\exp\Big\{
\sum_{j=1}^{\wh N(t)} \log E[
\e^{ixL(t-Z_j, t+s-Z_j]}
] 
\Big\}\mid \wh N(t) \Big] \\
& \qquad \times \exp\Big\{
-\int_{v=0}^1 \int_{r=t-v}^{t+s-v} (1-E[
\e^{ixL(t+s-v-r)}
])\Lambda(\d v) F_D(\d r)
\Big\} \\
& = (E[\e^{ixL(t-Z,t+s-Z)}])^{\wh N(t)} 
\exp\Big\{
-\int_{v=0}^1 \int_{r=t-v}^{t+s-v} (1-E[
\e^{ixL(t+s-v-r)}
]) \Lambda(\d v) F_D(\d r)
\Big\}, 
\end{align*}
where $(Z_j)$ is the iid random sequence with generic random element $Z$. 
The last expression coincides with the result and the proof is over. 
\end{proof}
\noindent
Note that due to the convolution $G:=\Lambda \ast F_D$, the last
term in \eqref{condi:ch:poisson:observed} has another expression
\begin{align*}
& \exp\Big\{
-\int_{v=0}^1 \int_0^\infty (1-E[
\e^{ixL(t+s-v-r)}1_{\{
t\le r+v \le t+s
\}}
]) \Lambda(\d v) F_D(\d r)
\Big\} \\
& = \exp\Big\{
-E[\wh N(t,t+s]] \int_t^{t+s}(1-E[\e^{ixL(t+s-w)}]) G(\d w)/ E[\wh N(t,t+s]]
\Big\} \\
& = E\Big[
\Big(E\big[
\e^{ixL(t+s-W)}
\big]
\Big)^{\wh N(t,t+s]}
\Big],
\end{align*}
where $W$ is independent of $L$ and has distribution $G(\d w)/ E[\wh
N(t,t+s]]$ on $(t,t+s]$. 

Now we differentiate the conditional ch.f. at $x=0$ proper
times to
obtain the following result. 

\bth\label{thm:condi:mamont:duration}
Consider the model \eqref{eq:cpp} with 
iid additive clusters 
$L_k,\,k=1,2,\ldots$ given by \eqref{chf:additive:process} such that cluster processes start at time
$(T_k+D_k)$. 
Let $\ell=0, 1,2,\ldots,\,t\ge 1$ and $s>0$.\\
$($i$)$ The prediction $\wt M_\ell(t,t+s]$ of $M(t,t+s]$ given $\{\wh
N(t)=\ell\}$ is given by 
\[
\wt M_\ell(t,t+s] = \ell J_1 +H_1,
\] 
where 
\[
J_i= 
\int_0^t E[L^i(t-u,t+s-u]] \frac{\Lambda \ast F_D(\d u)}{E[\wh N(t)]},\quad i=1,2
\]
and 
\[
 H_i=\int_0^1 \int_{r=t-v}^{t+s-v}E[L^i(t+s-v-r)]\Lambda(\d v)F_D(\d
 r),\quad i=1,2.
\]
$($ii$)$ The conditional variance of $M(t,t+s]$ given $\{\wh
N(t)=\ell\}$ is 
\begin{align*}
\mathrm{Var}(M(t,t+s]\mid \wh N(t)=\ell ) = \ell J_2 -\ell J_1^2 +H_2.
\end{align*}
\ethe

\bre
Since $E[\wh N(t)]=\Lambda(1)-\wh \Lambda(t)$, we evaluate the error of
prediction $\wt M$ by 
\[
E[(M(t,t+s]-E[M(t,t+s]\mid \wh N(t)])^2]= E[\mathrm{Var}(M(t,t+s]\mid \wh N(t))] = (\Lambda(1)-\wh
 \Lambda(t))(J_2-J_1^2) +H_2,
\]
which we could not do in the prediction by $\wh M$ of Section \ref{sec:pipcm}.
\ere

Applying Theorem \eqref{thm:condi:mamont:duration}, we calculate 
the following examples. 

\ble
Let $L$ be a NP process 
 with the 
mean value function $\mu(t)>0$. Then 
\[
\wt M_\ell (t,t+s] = \ell \int_0^t \mu(t-u,t+s-u]
\frac{\Lambda\ast F_D(\d u)}{E[\wh N(t)]} + \int_0^1
\int_{r=t-v}^{t+s-v}\mu(t+s-v-r)\Lambda(\d v) F_D(\d r)
\]
and the conditional variance of $M(t,t+s]$ given $\wh N(t)=\ell$ is 
\begin{align*}
\mathrm{Var}(M(t,t+s]\mid \wh N(t)=\ell ) &= \wt
M_\ell (t,t+s] -\ell \left(
\int_0^t \mu(t-u,t+s-u] \frac{\Lambda\ast F_D(\d u)}{E[\wh N(t)]} 
\right)^2\\
& + \ell 
\int_0^t \mu^2(t-u,t+s-u]  \frac{\Lambda\ast F_D(\d u)}{E[\wh N(t)]} \\
& + \int_0^1 \int_{r=t-v}^{t+s-v}\mu^2(t+s-v-r)\Lambda(\d v)F_D(\d r).
\end{align*}
\ele

\ble
Let $L$ be $NB$ process defined in \eqref{pb:NB}. Then the prediction $\wt
M_\ell (t,t+s]$ of $M(t,t+s]$ given $\{\wt N(t)=\ell \}$ is given by 
\[
\wt M_\ell (t,t+s] = \ell \frac{q}{p} \int_0^t \mu(t-u,t+s-u]
\frac{\Lambda\ast F_D(\d u)}{E[\wh N(t)]}+ \frac{q}{p}\int_0^1
\int_{r=t-v}^{t+s-v}\mu(t+s-v-r) \Lambda(\d v) F_D(\d r)
\]
and the conditional variance of $M(t,t+s]$ given $\wh N(t)=\ell$ is 
\begin{align*}
\mathrm{Var}(M(t,t+s]\mid \wh N(t)=\ell ) &= \frac{\wt
M_\ell (t,t+s] }{p} + \ell \frac{q^2}{p^2}\int_0^t \mu^2(t-u,t+s-u]
\frac{\Lambda\ast F_D(\d u)}{E[\wh N(t)]} \\
& -\ell \frac{q^2}{p^2} \left(\int_0^t \mu(t-u,t+s-u]
\frac{\Lambda\ast F_D(\d u)}{E[\wh N(t)]}\right)^2\\
& + \frac{q^2}{p^2} \int_0^1
\int_{r=t-v}^{t+s-v}\mu^2(t+s-v-r)\Lambda(\d v) F_D(\d r).
\end{align*}
\ele

\section{Numerical examples and some discussion}
\begin{figure}[htbp!]\label{figures}
\centerline{
\epsfig{figure=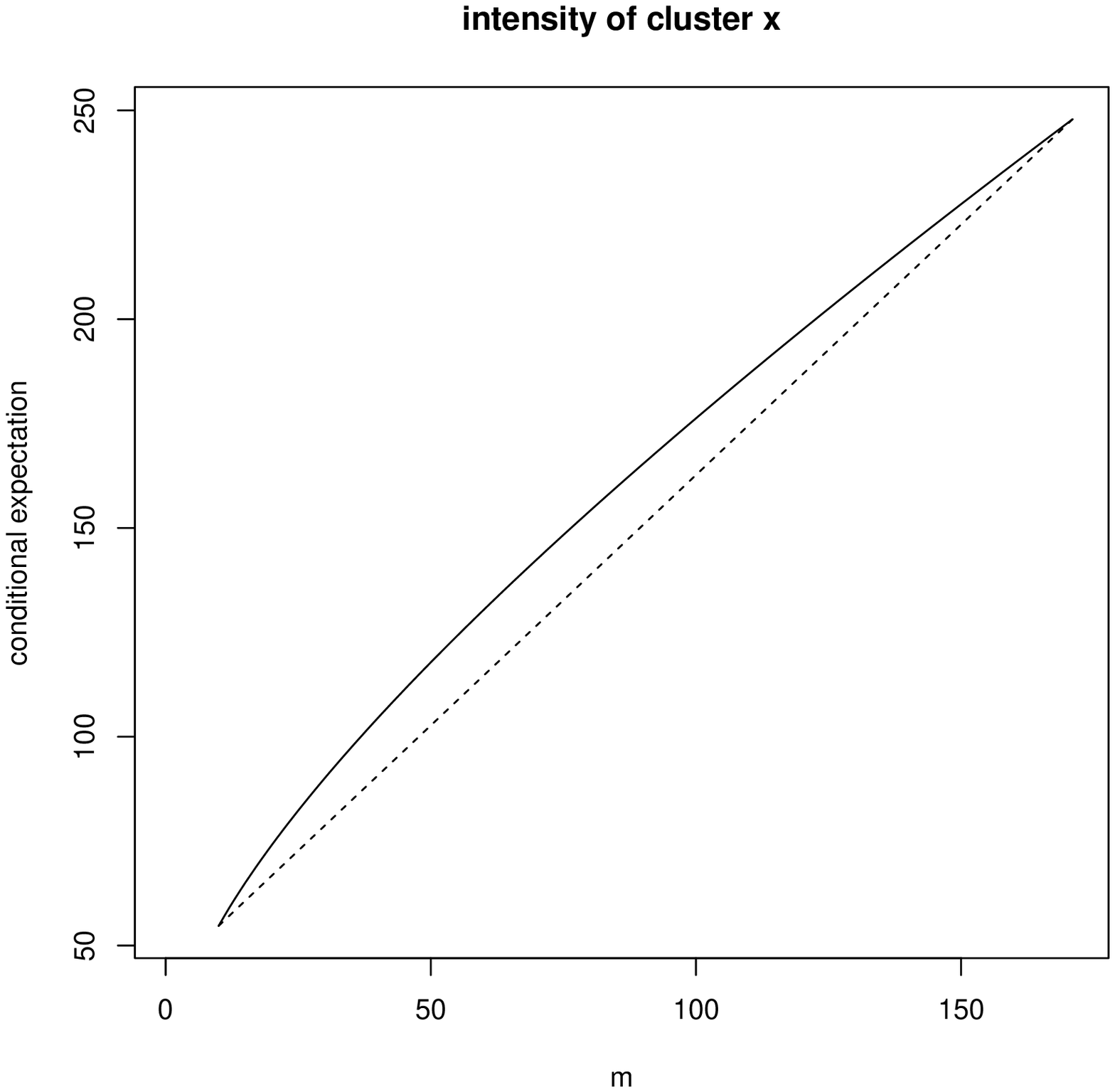, height=7cm, width=7cm}
\epsfig{figure=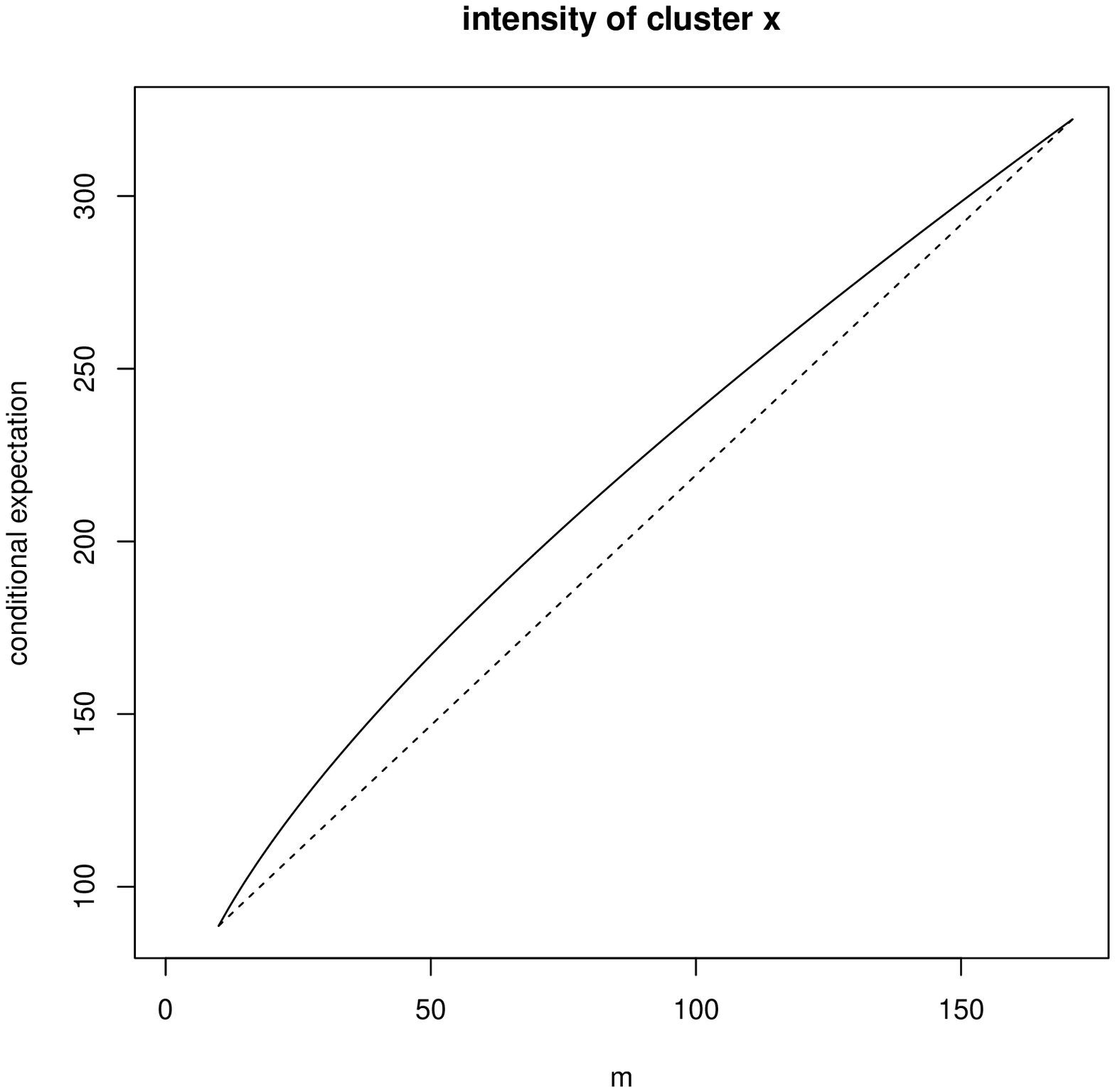,height=7cm,width=7cm}}
\centerline{
\epsfig{figure=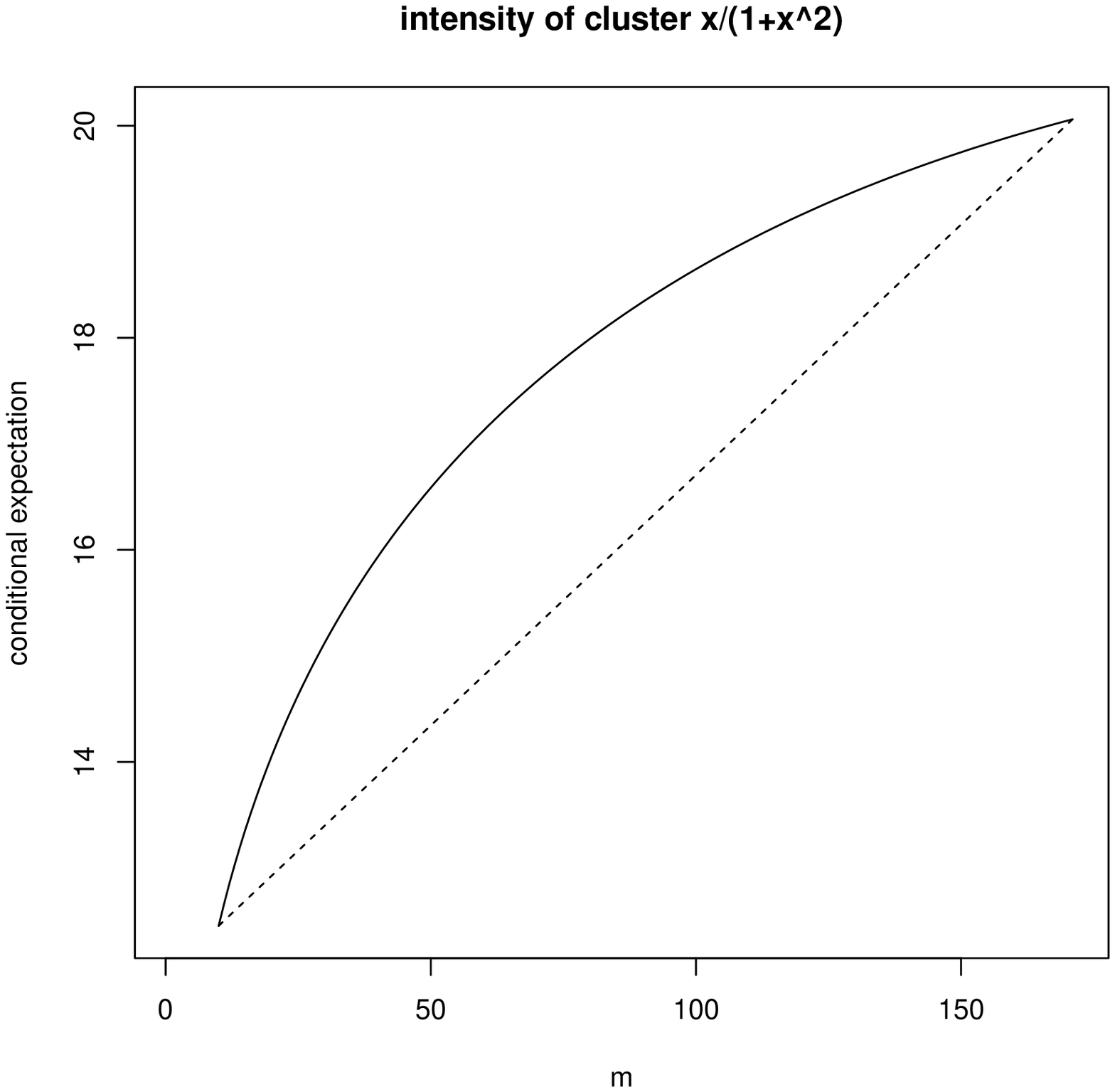,height=7cm,width=7cm}
\epsfig{figure=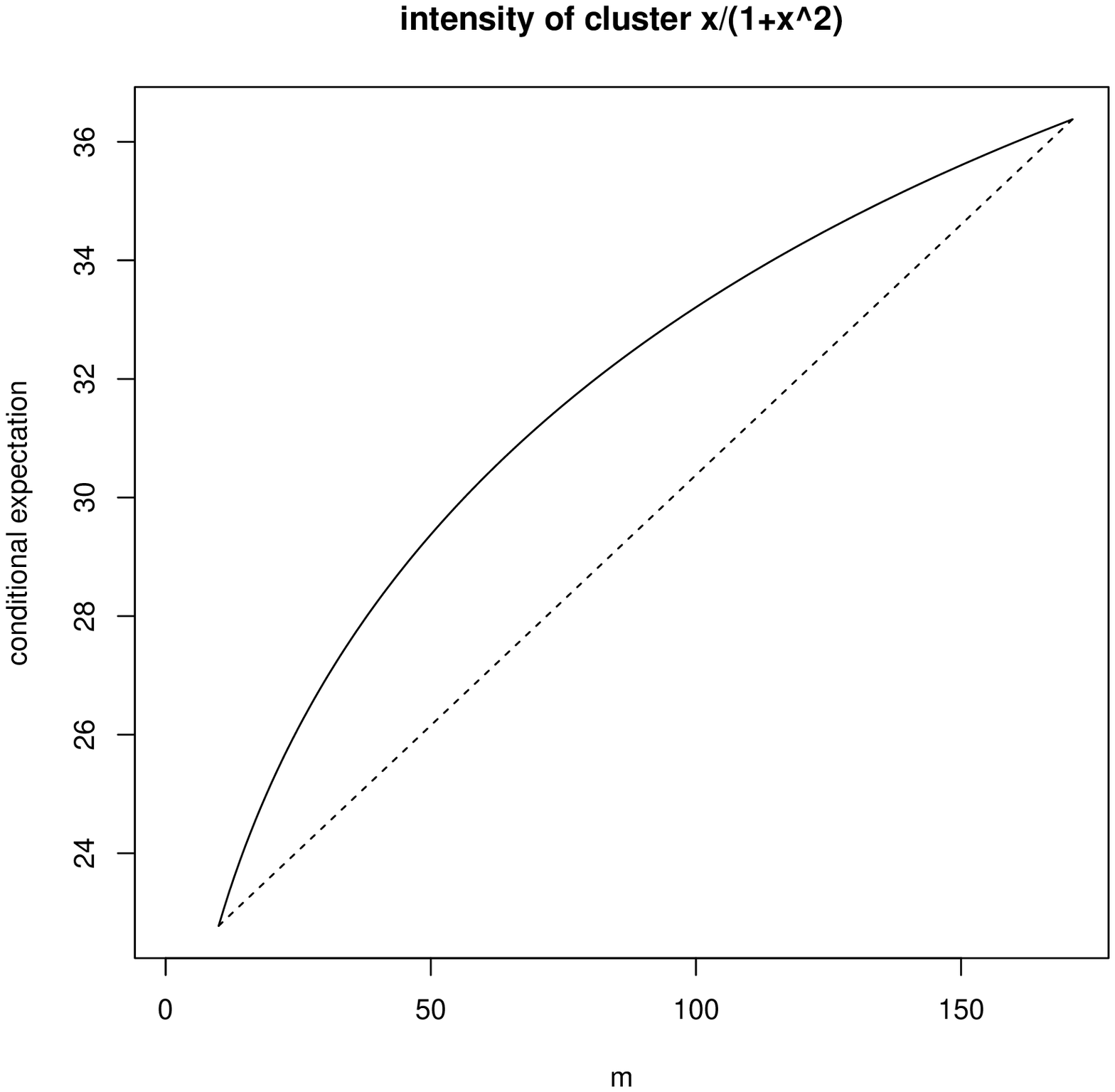,height=7cm,width=7cm}
}
\centerline{
\epsfig{figure=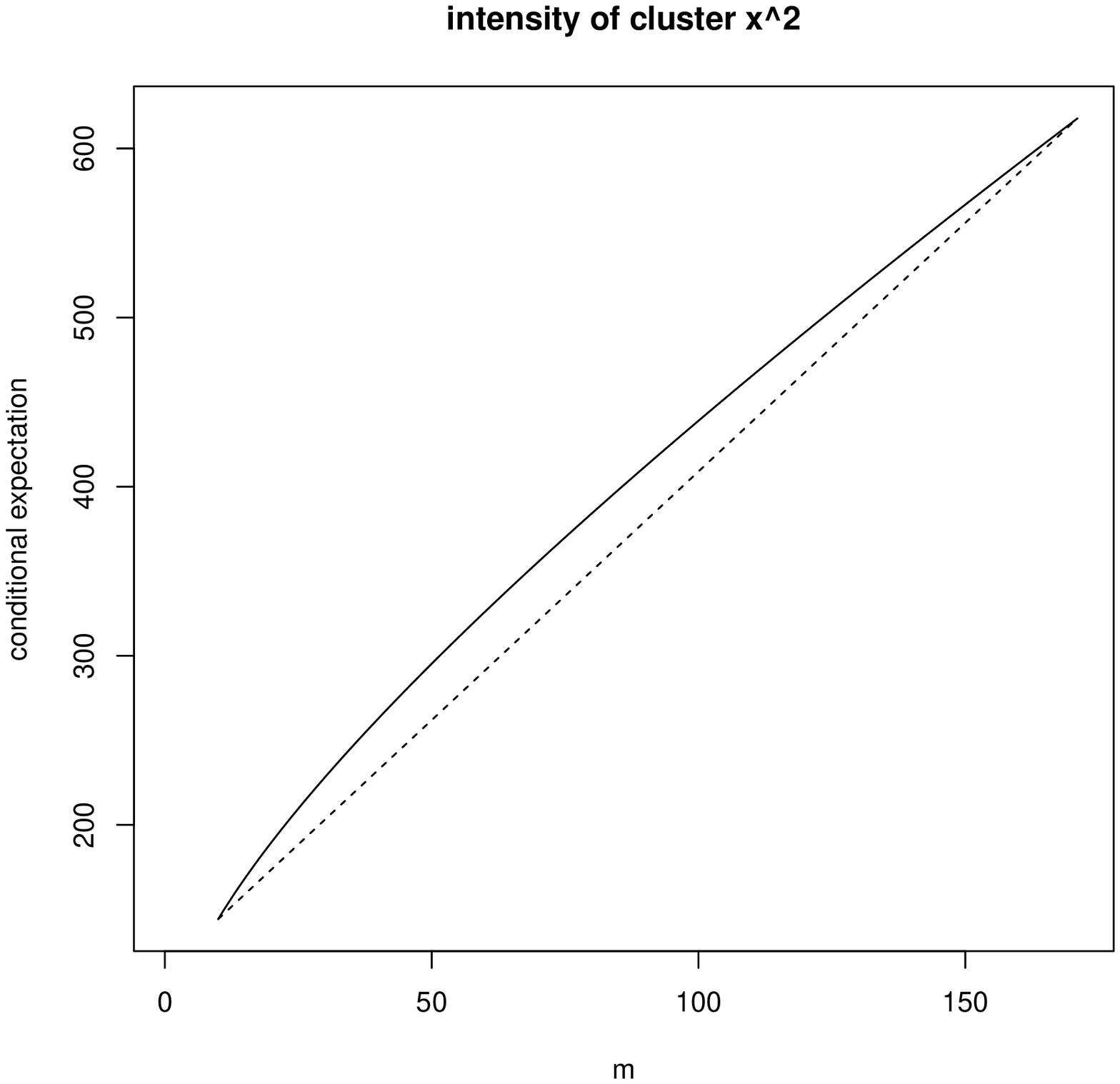,height=7cm,width=7cm}
\epsfig{figure=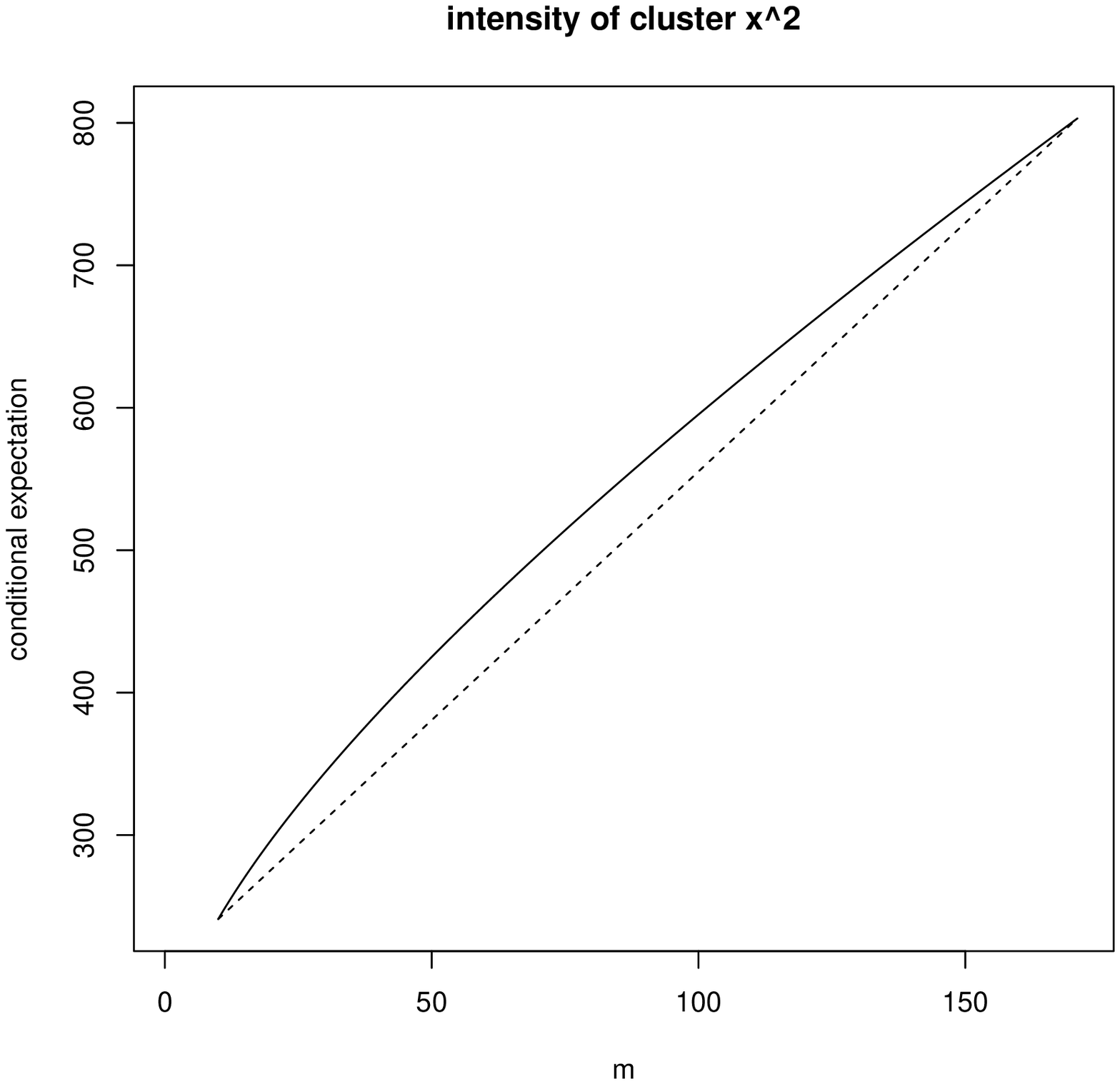,height=7cm,width=7cm}
}
\caption{
Graphs of conditional expectation of $M_m(1,2]$ based on the recursive
 algorithm given in Subsection \ref{subsec:nhpc}.
{\em Top left}, $(\Lambda_1,\mu_1)$. 
{\em Top right}, $(\Lambda_2,\mu_1)$.
{\em Middle left}, $(\Lambda_1,\mu_2)$. 
{\em Middle right}, $(\Lambda_2,\mu_2)$. 
{\em Bottom left}, $(\Lambda_1,\mu_3)$.
{\em Bottom right}, $(\Lambda_2,\mu_3)$.
}\label{fig:0}
\end{figure}

In this section we will observe how
non-homogeneity affects the predictor with several examples.
We consider the predictor $\widehat M_m(t,t+s]$ (see Subsection
\ref{subsec:nhpc}) with $NP$ 
clusters $L_j,\,j=1,2,\ldots$ under different mean value functions $\mu$
where we keep the underlying Poisson processes $N$ homogeneous. 
For the mean value function of the process $N$, two cases $E[N(x)]=\Lambda_1(x)=30x,\,\Lambda_2(x)=60x$ are examined, 
whereas we set three mean value functions for the cluster $L$, which are
\begin{align*}
 \mu_1(x)=5x, \qquad \mu_2(x)=\frac{5x}{1+x^2}, \quad \mathrm{and}\quad \mu_3(x)=5x^2.
\end{align*}
The middle one is a decreasing function while other two are
increasing ones. We plot the predictor $\widehat M_m(1,2]$ 
as function of $m$ for $m=10\sim 170$ in Figure \ref{figures}. We also make a
straight dot line from the initial value to the end value for comparison. 
In the light of Figure \ref{figures}, we see non-linearity of $\widehat
M_m(1,2]$ as a
function $m$ in all cases, and sizes of $\widehat M_m(1,2]$ properly
reflect the strength of intensity
functions.

Finally, we mention how our model \eqref{eq:cpp} could be estimated from data. 
The process $N(t)$ may be estimated from the claims arrivals observations,
whereas the generic process $L$ of clusters $(L_j)$ would be estimated
from observed payment streams. Nowadays statistical estimations of stochastic
processes are well established 
and since our model
uses basic processes which are not restrictive, we may have
no difficulty in estimation. Then once
the model \eqref{eq:cpp} is estimated, the prediction of future payment
amount would be possible by our proposed method. 

\end{document}